# EQUALITY OF CRITICAL POINTS FOR POLYMER DEPINNING TRANSITIONS WITH LOOP EXPONENT ONE

By Kenneth S. Alexander[1] and Nikos Zygouras

*University of Southern California and University of Warwick*

We consider a polymer with configuration modelled by the trajectory of a Markov chain, interacting with a potential of form $u + V_n$ when it visits a particular state 0 at time $n$, with $\{V_n\}$ representing i.i.d. quenched disorder. There is a critical value of $u$ above which the polymer is pinned by the potential. A particular case not covered in a number of previous studies is that of loop exponent one, in which the probability of an excursion of length $n$ takes the form $\varphi(n)/n$ for some slowly varying $\varphi$; this includes simple random walk in two dimensions. We show that in this case, at all temperatures, the critical values of $u$ in the quenched and annealed models are equal, in contrast to all other loop exponents, for which these critical values are known to differ, at least at low temperatures.

**1. Introduction.** A polymer pinning model is described by a Markov chain $(X_n)_{n\geq 0}$ on a state space containing a special point 0 where the polymer interacts with a potential. The space-time trajectory of the Markov chain represents the physical configuration of the polymer, with the $n$th monomer of the polymer chain located at $(n, X_n)$ (or just at $X_n$, for an undirected model). When the chain visits 0 at some time $n$, it encounters a potential of form $u + V_n$. The i.i.d. random variables $(V_n)_{n\geq 1}$ typically model variation in monomer species. We study the phase transition in which the polymer depins from the potential when $u$ goes below a critical value. We denote the distribution of the Markov chain (started from 0) in the absence of the potential by $P^X$ and we assume that it is recurrent. This recurrence assumption is merely a convenience and does not change the essential mathematics; see [1, 11]. Of greatest interest is the case where the excursion length

Received November 2008; revised June 2009.
[1]Supported by NSF Grants DMS-04-05915 and DMS-08-04934.

distribution decays as a power law:

$$P^X(\mathcal{E} = n) = n^{-c}\varphi(n), \qquad n \geq 1. \tag{1.1}$$

Here, the *loop exponent* is $c \geq 1$, $\mathcal{E}$ denotes the length of an excursion from 0, that is, the time elapsed between successive returns to 0, and $\varphi$ is a slowly varying function, that is, a function satisfying $\varphi(\kappa n)/\varphi(n) \to 1$ as $n$ tends to infinity for all $\kappa > 0$.

A large part of the existing rigorous literature on such models omits the case $c = 1$ because it is often technically different and not covered by the methods that apply to $c > 1$; see, for example, [1, 13, 15, 17]. That omission is partially remedied in this paper and we will see that the behavior for $c = 1$ can be quite different from the behavior for $c > 1$. The case $c = 1$ includes symmetric simple random walk in two dimensions, for which $\varphi(n) \sim \pi/(\log n)^2$ [14]. The essential feature of $c = 1$ is that $P^X(\mathcal{E} > n)$ is a slowly varying function of $n$ so that, for example, the longest of the first $m$ excursions typically has length greater than any power of $m$. This effectively enables the polymer to (at low cost) bypass stretches of disorder in which the values $V_n$ are insufficiently favorable and make returns to 0 in more favorable stretches.

The quenched version of the pinning model is described by the Gibbs measure

$$d\mu_N^{\beta,u,\mathbf{V}}(\mathbf{x}) = \frac{1}{Z_N} e^{\beta H_N^u(\mathbf{x},\mathbf{V})} dP^X(\mathbf{x}), \tag{1.2}$$

where $\mathbf{x} = (x_n)_{n \geq 0}$ is a path, $\mathbf{V} = (V_n)_{n \geq 0}$ is a realization of the disorder and

$$H_N^u(\mathbf{x}, \mathbf{V}) = \sum_{n=1}^{N} (u + V_n)\delta_0(x_n). \tag{1.3}$$

The normalization

$$Z_N = Z_N(\beta, u, \mathbf{V}) = E^X[e^{\beta H_N^u(\mathbf{x},\mathbf{V})}]$$

is the partition function. The disorder $\mathbf{V}$ is a sequence of i.i.d. random variables with mean zero and variance one. We denote the distribution of this sequence by $P^V$. We assume that $V_1$ has exponential moments of all orders and denote by $M_V(\beta)$ the moment generating function of $P^V$.

Let

$$L_N^X = L_N^X(\mathbf{x}) := \sum_{n=0}^{N} \delta_0(x_n)$$

denote the local time at 0 and define the *quenched free energy*

$$f_q(\beta, u) := \frac{1}{\beta} \lim_{N \to \infty} \frac{1}{N} \log Z_N(\beta, u, \mathbf{V}),$$

where this limit is taken $P^V$-a.s. The existence and nonrandomness of this limit is standard, as is the fact that

$$f_q(\beta, u) = \frac{1}{\beta} \lim_{N \to \infty} \frac{1}{N} E^V \log Z_N(\beta, u, \mathbf{V});$$

see [8]. The parameter $u \in \mathbb{R}$ can be thought of as the mean value of the potential, while the parameter $\beta > 0$ is the inverse temperature. It is known that the phase space in $(\beta, u)$ is divided by a critical line $u = u_c^q(\beta)$ into two regions: *localized* and *delocalized*. In the delocalized region $u < u_c^q(\beta)$, we have $f_q(\beta, u) = 0$, while in the localized region $u > u_c^q(\beta)$, we have $f_q(\beta, u) > 0$. It is proved in [12] that $f_q(\beta, \cdot)$ is infinitely differentiable for all $u > u_c^q(\beta)$. An alternate, more phenomenological, characterization of the two regions is as follows. From convexity, we have, for fixed $\beta$, that

(1.4) $$\left\langle \frac{L_N}{N} \right\rangle_N^{\beta, u, \mathbf{V}} = \frac{1}{\beta} \frac{\partial}{\partial u} \left( \frac{1}{N} \log Z_N(\beta, u, \mathbf{V}) \right) \to \frac{\partial}{\partial u} f_q(\beta, u) \qquad \text{for all } u,$$

$P^V$-a.s., where $\langle \cdot \rangle_N^{\beta, u, \mathbf{V}}$ denotes expectation under $\mu_N^{\beta, u, \mathbf{V}}$. This limiting value is called the *contact fraction*, denoted $C_q(\beta, u)$, and it is positive in the localized region and zero in the delocalized region. When the contact fraction is positive, we say the polymer is *pinned*.

The effect of the quenched disorder on the phase transition is quantified by comparing the quenched model to the corresponding *annealed model*, which is obtained by averaging the quenched Gibbs weight over the disorder to give the annealed Gibbs weight

$$E^V(e^{\beta H_N^u(\mathbf{x}, \mathbf{V})}) = e^{\beta \Delta L_N(\mathbf{x})},$$

where $\Delta = u + \beta^{-1} \log M_V(\beta)$. The corresponding annealed partition function is

$$Z_N^a = Z_N^a(\beta, u) := E^X(e^{\beta \Delta L_N})$$

and the Gibbs measure is

(1.5) $$d\mu_N^{\beta, u}(\mathbf{x}) = \frac{1}{Z_N^a} e^{\beta \Delta L_N(\mathbf{x})} dP^X(\mathbf{x}).$$

The corresponding annealed free energy and contact fraction are denoted $f_a(\beta, u)$ and $C_a(\beta, u)$, respectively. The annealed critical point is readily shown to be $u_c^a(\beta) = -\beta^{-1} \log M_V(\beta)$ for all $\beta > 0$ (see [2]), so $\Delta = u - u_c^a(\beta)$. It is a standard consequence of Jensen's inequality that $f_a(\beta, u) \geq f_q(\beta, u)$, so $u_c^a(\beta) \leq u_c^q(\beta)$. The effect, or lack of effect, of the disorder on the depinning transition may be seen in whether these two critical points actually differ and whether the specific heat exponent (describing the behavior of the free energy as $u$ decreases to the critical point) is different in the quenched case.

Although most mathematically rigorous work is relatively recent, there is an extensive physics literature on polymer pinning models; see the recent book [8] and the surveys [9, 16] and references therein. In [1] (see also [15] for a slightly weaker statement with simpler proof), it was proven that for $1 < c < 3/2$, and for $c = 3/2$ with $\sum_{n=1}^{\infty} 1/n\varphi(n)^2 < \infty$, for sufficiently small $\beta$, one has $u_c^q(\beta) = u_c^a(\beta)$ and the specific heat exponents are the same. Both works considered Gaussian disorder, although the method in [1] can be extended to accommodate more general disorder having a finite exponential moment.

By contrast, it follows straightforwardly from the sufficient condition ([17], (3.6)) that for $c > 1$, if $V_1$ is unbounded, or if $V_1$ is bounded and its essential supremum $v$ satisfies $P^V(V_1 = v) = 0$, then, for sufficiently large $\beta$, one has $u_c^q(\beta) > u_c^a(\beta)$; the method is based on fractional moment estimates. These results, together with [1], suggest that for $1 < c < 3/2$, there should be a transition from weak to strong disorder, that is, there should exist a value $\beta_0 > 0$ below which the annealed and quenched critical curves coincide [i.e., $u_c^q(\beta) = u_c^a(\beta)$ for $\beta < \beta_0$, while for $\beta > \beta_0$, one has $u_c^a(\beta) < u_c^q(\beta)$], but this has not been proven.

For $1 < c < 3/2$ and certain choices of bounded $V_1$ [necessarily with $P(V_1 = v) > 0$], it is known that the quenched and annealed critical points are equal for all $\beta > 0$ [5]. However, in these examples, $\operatorname{Var}(e^{\beta V_1})/[E(e^{\beta V_1})]^2$ stays bounded as $\beta \to \infty$, so there is no true "strong disorder" regime.

For $c > 3/2$, it follows from [11] that the quenched and annealed specific heat exponents are different and it was proven in [4] that the critical points are strictly different for all $\beta > 0$, that is, $\beta_0 = 0$. In [3], the distinctness of critical points at high temperature was extended to include $c = 3/2$ with $\varphi(n) \to 0$ as $n \to \infty$ and the asymptotic order of the gap $u_c^q(\beta) - u_c^a(\beta)$ was given. Recently, in [10], the critical points were shown to be distinct for all $\beta > 0$ for the case of $c = 3/2$ and $\varphi(n)$ asymptotically a positive constant, a case about which physicists had long disagreed [6, 7].

Here, we show that even with true strong disorder, the critical points remain the same in the case $c = 1$.

THEOREM 1.1. *Consider the quenched model (1.2) and suppose that $E(e^{tV_1}) < \infty$ for all $t \in \mathbb{R}$ and that (1.1) holds with $c = 1$. For all $\beta > 0$ and all $u > u_c^a(\beta)$, the quenched free energy $f_q(\beta, u) > 0$ and thus $u_c^q(\beta) = u_c^a(\beta)$ for all $\beta > 0$.*

In [1] and [15], for the case $1 < c < 3/2$, a statement stronger than the equality of the critical points was proven: given $\varepsilon > 0$, if $\beta$ and $\beta\Delta$ are sufficiently small, then one has $f_a(\beta, u) \geq f_q(\beta, u) > (1-\varepsilon)f_a(\beta, u)$. One may ask whether a similar statement (possibly strengthened to be valid for all $\beta$'s) holds for the $c = 1$ case. We do not pursue that question here, although we expect such a statement to be true for small $\beta$. There are technical obstacles to carrying over the proof for $1 < c < 3/2$ to the case $c = 1$, as noted in Section 4 of [1].

**2. Notation and idea of the proof.** Denote the local time at zero over a time interval $I$ by

$$L_I^X = \sum_{n \in I} \delta_0(x_n), \tag{2.1}$$

so that $L_N^X = L_{[0,N]}^X$. The overlap between two paths $X, X'$ in an interval $I$ is defined as

$$B_I^{X,X'} = \sum_{n \in I} \delta_0(x_n)\delta_0(x_n'). \tag{2.2}$$

We denote by $P^{X,X'}$ the measure corresponding to two independent copies $X, X'$ of the Markov chain. The "energy gained over an interval $I$" is defined as

$$H_I^u(\mathbf{x}, V) = \sum_{n \in I}(u + V_n)\delta_0(x_n). \tag{2.3}$$

The annealed correlation length is defined to be $M = M(\beta, u) := 1/(\beta f_a(\beta, u))$. From (1.5), both $\beta f_a(\beta, u)$ and $M$ are functions of only the product $\beta\Delta$. Using Laplace asymptotics and the large deviations for the local time $L_N$, one can deduce the asymptotics of $M$ and $C_a(\beta, u)$ for $\beta\Delta \to 0$. Specifically, letting

$$\Psi(t) = \int_t^\infty \varphi(e^s)\, ds,$$

we obtain

$$\beta\Delta \sim \Psi(\log M) \quad \text{and} \quad C_a(\beta, u) \sim \frac{1}{M\varphi(M)} \qquad \text{as } \beta\Delta \to 0.$$

For example, if $\varphi(n) \sim K(\log n)^{-\alpha}$ for some $\alpha > 1$, then

$$\log M = \log \frac{1}{\beta f_a(\beta, u)} \sim \left(\frac{\alpha - 1}{K}\beta\Delta\right)^{-1/(\alpha-1)} \qquad \text{as } \beta\Delta \to 0, \tag{2.4}$$

so $f_a(\beta, \cdot)$ is $C^\infty$, even at $u = u_c^a(\beta)$. The details are similar to those in the case $c > 1$ considered in [1], but we do not include them here as they are not required for our analysis.

We use length scales $K_1(\beta, M), K_2(\beta, M)$, related as follows, for $\beta, M > 0$. Let $\Lambda_V(\beta) := \log M_V(2\beta) - 2\log M_V(\beta)$. For $c = 1$, (1.1) implies that $\varphi(x) \to 0$ as $x \to \infty$. Since $\varphi$ is slowly varying, this in turn implies that

$$\frac{\log x}{\log(1/\varphi(x))} \to \infty \qquad \text{as } x \to \infty.$$

Therefore, we can choose $K_1, K_2$ satisfying

(2.5) $$32K_2 < e^{\Lambda_V(\beta)K_2}$$

and

(2.6) $$4(M \vee 1)\log \frac{1}{\varphi(K_1)} < K_2 < \frac{1}{2\Lambda_V(\beta)}\log \frac{K_1}{2}.$$

For fixed $\beta$, as $\Delta \to 0$ (i.e., $M \to \infty$), we then have $M \ll K_2 \ll K_1$. We assume henceforth that $K_1, K_2$ are even integers.

Define the intervals
$$I_i = [iK_1, (i+1)K_1) \cap \mathbb{Z}, \qquad I_i^\gamma = [iK_1, (i+\gamma)K_1) \cap \mathbb{Z}$$
for $0 < \gamma < 1$. For an interval $I$, let $\tau_I = \inf\{n \in I : x_n = 0\}$ and $\sigma_I = \sup\{n \in I : x_n = 0\}$. We set $\tau_I = \sigma_I = \infty$ if the path does not visit 0 during the interval $I$. We denote by $\Xi_{NK_1}$ the set of all paths of length $NK_1$ which have the following property: if $\tau_{I_i} < \infty$ for some $i \leq N$, then $\tau_{I_i} \in I_i^{1/2}$ and $\sigma_{I_i} - \tau_{I_i} \leq K_2$.

*Idea of the proof.* We will look at a scale $NK_1$ and restrict the partition function $Z_{NK_1}(u, \beta, \mathbf{V})$ to paths that belong to the set $\Xi_{NK_1}$. Further, we will restrict our attention to paths within $\Xi_{NK_1}$ which bypass *bad* blocks of length $K_1$. Roughly speaking, a bad block is defined to be a block for which the quenched partition function of a path starting at a uniform random point in the block, and making its final visit to 0 in the block within time $K_2$ after this starting point, is less than half of the corresponding annealed partition function. In Lemma 3.2, we control the probability of having a bad block. It then remains to make an energy-entropy balancing of the paths that belong in $\Xi_{NK_1}$ and bypass bad blocks, and to show that for $\beta > 0$ and $\Delta = u + \beta^{-1}\log M_V(\beta) > 0$, this balance is uniformly (in $N$) bounded away from zero. For this, we will use the fact that in a *good* block, the free energy gained is of the order $K_2/M$ (this is essentially Lemma 3.1), and the fact that because $P^X(\mathcal{E} > k)$ is a slowly varying function of $k$, the cost of bypassing bad blocks is small.

## 3. Proof of the theorem.

LEMMA 3.1. *Let $\beta > 0, u \in \mathbb{R}, \Delta = u + \beta^{-1}\log M_V(\beta)$ and $M = M(\beta, u)$. Then, for all $N > 2\beta\Delta M$,*

(3.1) $$\log E^X[e^{\beta\Delta L_N}] \geq \frac{1}{2}\frac{N}{M}.$$

PROOF. It is observed in [1] that $a_N := \beta\Delta + \log E^X[e^{\beta\Delta L_N}]$ is subadditive in $N$. Since $a_N/N \to \beta f_a(\beta, u)$, it follows that
$$\beta\Delta + \log E^X[e^{\beta\Delta L_N}] \geq N\beta f_a(\beta, \Delta) = \frac{N}{M}$$

and the result is immediate. □

The block $I_i$ is called *good* if it satisfies

$$\sum_{b \in I_i^{1/2}} E^X[e^{\beta H^u_{[b,b+K_2]}(\mathbf{x},\mathbf{V})}|x_b=0] > \frac{1}{2} \sum_{b \in I_i^{1/2}} E^V E^X[e^{\beta H^u_{[b,b+K_2]}(\mathbf{x},\mathbf{V})}|x_b=0]$$

$$= \frac{|I_i^{1/2}|}{2} E^X[e^{\beta \Delta L_{K_2}}]$$

and called *bad* otherwise. Let $p^V_{\text{good}} := P^V(I_i \text{ is good})$ and $p^V_{\text{bad}} := P^V(I_i \text{ is bad})$.

PROPOSITION 3.2. *For $K_1, K_2$ satisfying (2.5) and the second inequality in (2.6), we have $p^V_{\text{good}} > 1/2$.*

PROOF. By Chebyshev's inequality,

$$p^V_{\text{bad}} \le 4 \frac{\text{Var}^V(\sum_{b \in I_i^{1/2}} E^X[e^{\beta H^u_{[b,b+K_2]}(\mathbf{x},\mathbf{V})}|x_b=0])}{(\sum_{b \in I_i^{1/2}} E^V E^X[e^{\beta H^u_{[b,b+K_2]}(\mathbf{x},\mathbf{V})}|x_b=0])^2}$$

$$< 4\Bigg( \sum_{b,b' \in I_i^{1/2}} 1_{|b-b'| \le K_2}$$

$$\times E^V E^{X,X'}[e^{\beta H^u_{[b,b+K_2]}(\mathbf{x},\mathbf{V})+\beta H^u_{[b',b'+K_2]}(\mathbf{x}',\mathbf{V})}|x_b=x'_{b'}=0]\Bigg)$$

$$\times \Bigg( \sum_{b,b' \in I_i^{1/2}} E^{V,V'} E^{X,X'}[e^{\beta H^u_{[b,b+K_2]}(\mathbf{x},\mathbf{V})(X)+\beta H^u_{[b',b'+K_2]}(\mathbf{x}',\mathbf{V}')}|$$

$$x_b = x'_{b'} = 0]\Bigg)^{-1}.$$

Here, we used the fact that whenever the two independent paths $\mathbf{x}, \mathbf{x}'$ visit zero at points $b, b'$ such that $|b - b'| > K_2$, the energies $H^u_{[b,b+K_2]}(\mathbf{x},\mathbf{V})$ and $H^u_{[b',b'+K_2]}(\mathbf{x}',\mathbf{V})$ are independent.

An easy calculation shows that the above is equal to

$$4\Bigg( \sum_{b,b' \in I_i^{1/2}} 1_{|b-b'| \le K_2} E^{X,X'}[e^{\beta \Delta(L^X_{[b,b+K_2]}+L^{X'}_{[b',b'+K_2]})} e^{\Lambda_V(\beta) B^{X,X'}_{[b,b+K_2] \cap [b',b'+K_2]}}|$$

$$x_b = x'_{b'} = 0]\Bigg)$$

$$\times \left( \sum_{b,b' \in I_i^{1/2}} E^{X,X'}[e^{\beta \Delta(L^X_{[b,b+K_2]} + L^{X'}_{[b',b'+K_2]})} | x_b = x'_{b'} = 0] \right)^{-1}$$

$$\leq 4 \Bigg( \sum_{b,b' \in I_i^{1/2}} 1_{|b-b'| \leq K_2} e^{\Lambda_V(\beta) K_2}$$

$$\times E^{X,X'}[e^{\beta \Delta(L^X_{[b,b+K_2]} + L^{X'}_{[b',b'+K_2]})} | x_b = x'_{b'} = 0] \Bigg)$$

$$\times \left( \sum_{b,b' \in I_i^{1/2}} E^{X,X'}[e^{\beta \Delta(L^X_{[b,b+K_2]} + L^{X'}_{[b',b'+K_2]})} | x_b = x'_{b'} = 0] \right)^{-1}$$

$$= \frac{4}{|I_i^{1/2}|^2} \sum_{b,b' \in I_i^{1/2}} 1_{|b-b'| \leq K_2} e^{\Lambda_V(\beta) K_2} < \frac{32 K_2}{K_1} e^{\Lambda_V(\beta) K_2}$$

$$< \frac{1}{K_1} e^{2\Lambda_V(\beta) K_2} < \frac{1}{2}$$

for $K_1, K_2$ satisfying (2.5) and the second inequality in (2.6). In the third line, we have used the fact that the expectations in the second line do not depend on $b$ and $b'$. □

We now return to the proof of Theorem 1.1. Let

$$J_N := \{i \leq N : I_i \text{ is good}\} \cup \{0\} = \{i_1 < \cdots < i_{|J_N|}\}.$$

Under $P^V$, the sequence $(i_j - i_{j-1})_{j \geq 1}$ is an i.i.d. sequence of geometric random variables with parameter $p^V_{\text{good}}$.

We denote by $\Xi^{J_N}_{NK_1} = \Xi^{J_N}_{NK_1}(\mathbf{V})$ the set of paths $\mathbf{x} \in \Xi_{NK_1}$ which satisfy $x_{NK_1} = 0$ and make no returns to $0$ in bad blocks after the first block. In the following computation, $a_j$ and $b_j$ are the starting and ending points, respectively, of the excursion from $I_{i_j}$ to $I_{i_{j+1}}$. Let $p_n = P^X(\mathcal{E} = n)$. As a convention, we set $b_0 := 0$ and $b_{|J_N|} := NK_1$. Let $Z_{NK_1}(\Xi^{J_N}_{NK_1})$ denote the partition function restricted to the set of paths $\Xi^{J_N}_{NK_1}$. We then have

$$Z_{NK_1}(\Xi^{J_N}_{NK_1})$$
$$= \sum_{a_1 \leq K_2} \sum_{b_1 \in I_{i_2}^{1/2}} \sum_{a_2 - b_1 \leq K_2} \cdots \sum_{b_{|J_N|-1} \in I_{i_{|J_N|}}^{1/2}} \sum_{a_{|J_N|} - b_{|J_N|-1} \leq K_2}$$

$$\prod_{j=1}^{|J_N|} E^X[e^{\beta H^u_{[b_{j-1},a_j]}(\mathbf{x},\mathbf{V})} \delta_0(x_{a_j}) | x_{b_{j-1}} = 0] p_{b_j - a_j}$$

$$\geq \sum_{a_1 \leq K_2} \sum_{b_1 \in I_{i_2}^{1/2}} \sum_{a_2-b_1 \leq K_2} \cdots \sum_{b_{|J_N|-1} \in I_{i_{|J_N|}}^{1/2}} \sum_{a_{|J_N|}-b_{|J_N|-1} \leq K_2}$$

$$\prod_{j=1}^{|J_N|} E^X[e^{\beta H_{[b_{j-1},a_j]}^u(\mathbf{x},\mathbf{V})}; \sigma_{[b_{j-1},b_{j-1}+K_2]} = a_j | x_{b_{j-1}} = 0] p_{b_j-a_j}.$$

With a mild abuse of notation, let us interpret $I_{i_{|J_N|+1}}$ as meaning the one-point interval $\{NK_1\}$. On the set $\{\sigma_{[b_{j-1},b_{j-1}+K_2]} = a_j\}$, we have that $H_{[b_{j-1},a_j]}^u(\mathbf{x},\mathbf{V}) = H_{[b_{j-1},b_{j-1}+K_2]}^u(\mathbf{x},\mathbf{V})$ and, therefore, for some $C$, the above is bounded below by

$$\sum_{a_1 \leq K_2} \sum_{b_1 \in I_{i_2}^{1/2}} \sum_{a_2-b_1 \leq K_2} \cdots \sum_{b_{|J_N|-1} \in I_{|J_N|}^{1/2}} \sum_{a_{|J_N|}-b_{|J_N|-1} \leq K_2}$$

$$\prod_{j=1}^{|J_N|} E^X[e^{\beta H_{[b_{j-1},b_{j-1}+K_2]}^u(\mathbf{x},\mathbf{V})}; \sigma_{[b_{j-1},b_{j-1}+K_2]} = a_j | x_{b_{j-1}} = 0]$$

$$\times \min_{a \in I_{i_j}^{3/4}, b \in I_{i_{j+1}}^{1/2}} p_{b-a}$$

$$= E^X[e^{\beta H_{[0,K_2]}^u(\mathbf{x},\mathbf{V})}] \left( \min_{a \in I_{i_1}^{3/4}, b \in I_{i_2}^{1/2}} p_{b-a} \right)$$

$$\times \prod_{j=2}^{|J_N|} \sum_{b_{j-1} \in I_{i_j}^{1/2}} E^X[e^{\beta H_{[b_{j-1},b_{j-1}+K_2]}^u(\mathbf{x},\mathbf{V})} | x_{b_{j-1}} = 0] \min_{a \in I_{i_j}^{3/4}, b \in I_{i_{j+1}}^{1/2}} p_{b-a}$$

$$\geq E^X[e^{\beta H_{[0,K_2]}^u(\mathbf{x},\mathbf{V})}] \left( \min_{a \in I_{i_1}^{3/4}, b \in I_{i_2}^{1/2}} p_{b-a} \right)$$

$$\times \prod_{j=2}^{|J_N|} |I_{i_j}^{1/2}| E^X[e^{\beta \Delta L_{K_2}}] \min_{a \in I_{i_j}^{3/4}, b \in I_{i_{j+1}}^{1/2}} p_{b-a}$$

$$\geq E^X[e^{\beta H_{[0,K_2]}^u(\mathbf{x},\mathbf{V})}] C \frac{\varphi((i_2-i_1+1)K_1)}{(i_2-i_1+1)K_1}$$

$$\times \prod_{j=2}^{|J_N|} \left( C \frac{\varphi((i_{j+1}-i_j+1)K_1)}{(i_{j+1}-i_j+1)K_1} |I_{i_j}^{1/2}| E^X[e^{\beta \Delta L_{K_2}}] \right)$$

$$= \frac{1}{K_1} E^X[e^{\beta H_{[0,K_2]}^u(\mathbf{x},\mathbf{V})}]$$

$$\times \left[\prod_{j=1}^{|J_N|} C \frac{\varphi((i_{j+1} - i_j + 1)K_1)}{4(i_{j+1} - i_j + 1)}\right] (E^X[e^{\beta \Delta L_{K_2}}])^{|J_N|-1}.$$

In the second inequality, we used the fact that the interval $I_{i_j}$ is good, while the last equality makes essential use of $c = 1$ in the cancellation of factors $K_1$. We then have that

$$\frac{1}{NK_1} \log Z_{NK_1} \geq \frac{1}{NK_1} \log Z_{NK_1}(\Xi_{NK_1}^{J_N})$$

$$\geq \frac{1}{NK_1} \log \left(\frac{1}{K_1} E^X[e^{\beta H_{[0,K_2]}^u(\mathbf{x},\mathbf{V})}]\right) + \frac{|J_N|-1}{NK_1} \log E^X[e^{\beta \Delta L_{K_2}}]$$

$$+ \frac{1}{NK_1} \sum_{j=1}^{|J_N|} \log \frac{C\varphi((i_{j+1} - i_j + 1)K_1)}{4(i_{j+1} - i_j + 1)}.$$

Letting $N \to \infty$, we get that the left-hand side converges to the quenched free energy $f_q(\beta, u)$, while the right-hand side converges to

$$\frac{1}{K_1} p_{\text{good}}^V \log E^X[e^{\beta \Delta L_{K_2}}] + \frac{1}{K_1} p_{\text{good}}^V E^V \log \frac{C\varphi(i_2 K_1)}{i_2},$$

where $C$ is a constant different from what appears above. Recall that $i_2 - 1$ is a geometric random variable under $P^V$ with parameter $p_{\text{good}}^V$. For $K$ sufficiently large, we have

$$C_\varphi := \inf\left\{\frac{x\varphi(kx)}{\varphi(k)} : x \geq 1, k \geq K\right\} > 0$$

and we may assume that $K_1 \geq K$. We then have

$$f_q(\beta, u) \geq \frac{1}{K_1} p_{\text{good}}^V \left(\log E^X[e^{\beta \Delta L_{K_2}}] + E^V \log \frac{CC_\varphi \varphi(K_1)}{i_2^2}\right)$$

$$= \frac{1}{K_1} p_{\text{good}}^V (\log E^X[e^{\beta \Delta L_{K_2}}] + \log(CC_\varphi \varphi(K_1)) - 2E^V[\log i_2])$$

$$\geq \frac{1}{K_1} p_{\text{good}}^V \left(\log E^X[e^{\beta \Delta L_{K_2}}] + \log(CC_\varphi \varphi(K_1)) - 2\log\left(\frac{1}{p_{\text{good}}^V} + 1\right)\right)$$

and, by Lemma 3.1 and Proposition 3.2, this is bounded below by

$$(3.2) \qquad \frac{1}{2K_1}\left(\frac{K_2}{2M} + \log(CC_\varphi \varphi(K_1)) - 2\log 3\right).$$

Then, using the first inequality in (2.6), we get that, provided $M$ is sufficiently large, that is, $\Delta$ is small,

$$(3.3) \qquad f_q(\beta, u) > \frac{1}{2K_1}\left(\frac{K_2}{4M} + \log \frac{CC_\varphi}{9}\right) > 0.$$

This completes the proof of Theorem 1.1.

In the case $\varphi(n) \sim K(\log n)^{-\alpha}$ for some $\alpha > 1$, by (2.4), there are constants $C_i(\beta)$ such that for $\Delta < C_1(\beta)$, (2.5) and (2.6) are satisfied for $K_1 = \exp(C_2(\beta)M \log M)$ and $K_2 = C_3(\beta)M \log M$. Thus, the lower bound (3.3) says that

$$f_q(\beta, u) > e^{-C_4(\beta)M \log M}$$

with $M$ given approximately by (2.4). We know of no reason to believe that this bound is sharp.

DEPARTMENT OF MATHEMATICS KAP 108
UNIVERSITY OF SOUTHERN CALIFORNIA
LOS ANGELES, CALIFORNIA 90089-2532
USA
E-MAIL: alexandr@usc.edu

DEPARTMENT OF STATISTICS
UNIVERSITY OF WARWICK
COVENTRY CV4 7AL
UNITED KINGDOM
E-MAIL: N.Zygouras@warwick.ac.uk